\documentclass[11pt]{amsart}

\usepackage {amsfonts}
\usepackage[english]{babel}
\usepackage[cp1251]{inputenc}

\usepackage{amsmath}
\usepackage{amssymb}

\def\g{\gamma}
\def\G{\Gamma}
\def\d{\delta}
\def\a{\alpha}
\def\b{\beta}
\def\p{\varphi}
\def\e{\varepsilon}
\def\l{\lambda}
\def\L{\Lambda}
\def\o{\omega}
\def\O{\Omega}
\def\t{\theta}

\def\la{\langle}
\def\ra{\rangle}

\def\D{\mathcal D}
\def\S{\mathcal S}
\def\cN{\mathcal N}
\def\B{~\hfill$\Box$}
\def\R{{\mathbb R}}
\def\C{{\mathbb C}}
\def\N{{\mathbb N}}
\def\Z{{\mathbb Z}}

\DeclareMathOperator{\supp}{supp}

\newtheorem{Def}{Definition}
\newtheorem{Th}{Theorem}
\newtheorem{Pro}{Proposition}
\newtheorem*{Cor}{Corollary}

\begin{document}

\title{Growth of masses of crystalline measures}

\author{Peter Boyvalenkov,  Sergii Yu.Favorov}

\address{Peter Boyvalenkov,
\newline\hphantom{iii}  Institute of Mathematics and Informatics, Bulgarian Academy of Sciences,
\newline\hphantom{iii} 8 Acad. G. Bonchev str., Sofia 1113, Bulgaria}
\email{peter@math.bas.bg}

\address{Sergii Favorov,
\newline\hphantom{iii}  V.N.Karazin Kharkiv National University
\newline\hphantom{iii} Svobody sq., 4, Kharkiv, Ukraine 61022}
\email{sfavorov@gmail.com}

\maketitle {\small
\begin{quote}
\centerline {\bf Abstract.}

Let $\mu$ be a measure on the Euclidean space $\R^d$ of unbounded total variation that is positive or translation bounded and has
a pure point Fourier transform in the sense of distributions $\hat\mu$. We prove that the measure $\nu$ with the same support as $\hat\mu$ and masses equal to the squares of the masses of $\hat\mu$ is translation bounded.
We also prove that if $\mu$ is as above and the restriction of its spectrum, i.e., of the support of $\hat\mu$, to each ball of fixed radius is a linearly independent set over $\Z$, then the measure $\hat\mu$ is also translation bounded. These results imply certain conditions for a crystalline measure to be a Fourier quasicrystal.


\end{quote}
}

   \bigskip
\section{Introduction}\label{S1}
   \bigskip

 A complex measure $\mu$ on the Euclidean space $\R^d$ with a locally finite support (that is, its intersection with any compact set is finite)
 is called {\it crystalline}   if $\mu$ is a tempered distribution and its Fourier transform in the sense of distributions $\hat\mu$
 is also a measure with locally finite support. If both measures $|\mu|$ and $|\hat\mu|$ are tempered distributions,
then $\mu$ is called a {\it Fourier quasicrystal}. Here and below we denote by $|\nu|(E)$ the variation of the complex measure $\nu$ on the set $E$.

Recently, Fourier quasicrystals and crystalline measures are studied very actively. Many works are devoted to the investigations of their properties (see, for example, the survey papers \cite{M4,LT,AKV}). In fact,  Fourier quasicrystals are the form of Poisson formulas (see below), the latter were used in particular by D. Radchenko and M. Viazovska in \cite{RV}.

Special attention was paid to Fourier quasicrystals of the form
\begin{equation}\label{a}
\mu=\sum_{\l\in\L}\d_\l,\quad \l\in\R^d,
\end{equation}
 ($\d_\l$, as usual, means the unit mass at the point $\l$). A. Olevskii and A. Ulanovskii in \cite{OU}, \cite{OU1} for the case $d=1$ established 1-1 connection between such measures
and sets of zeros of real-rooted exponential polynomials. W. Lawton and A. Tsikh \cite{LT} partially extended this result in an arbitrary dimension $d$.
 In \cite{F6} the result of Olevskii and Ulanovskii was generalized to Dirichlet series. F. Goncalves \cite{G} investigated arbitrary measures $\mu$ on $\R$  with pure point measures $\hat\mu$ and obtained an analog of the result of Olevskii and Ulanovskii. Measures of the form \eqref{a} with $\l$ in a horizontal strip of finite width were studied in \cite{F7}.

Returning to the definition of a Fourier quasicrystal, it is natural to ask whether we need additional conditions on the variations
of the measures $\mu$ and $\hat\mu$ in order to prove that a crystalline measure is a Fourier quasicrystal.  The answer is positive, as in \cite{F5}
the second author constructed a crystalline measure which is not a Fourier quasicrystal.

In this paper, we investigate various additional conditions under which a crystalline measure becomes a Fourier quasicrystal. To this end, we investigate the properties of arbitrary tempered measures in $\R^d$ with a purely point Fourier transform $\hat\mu$. For $d=1$, this is exactly the class of measures that Goncalves considered in \cite{G}, but our investigation does not have other common points with \cite{G}.

To formulate our results, we recall some notions and properties related to the Fourier transform (see, e.g., \cite{Ru}):

Denote by $\D$ the space of $C^\infty$-functions on $\R^d$ with compact support and by $\S$ the Schwartz space of
 $C^\infty$-functions $\p(t)$ on $\R^d$ with finite norms
$$
\cN_m(\p)=\sup_{\R^d}\{\max\{1,|t|^m \} \cdot \max_{\|k\|\le m} |D^k\p(t)|\},\quad m=0,1,2,\dots,
 $$
where
$$
k=(k_1,\dots,k_d)\in(\N\cup\{0\})^d,\ \ \|k\|=k_1+\dots+k_d,\  \ D^k=\partial^{k_1}_{x_1}\dots\partial^{k_d}_{x_d}.
$$
The Fourier transform
$$
   \hat\p(y)=\int_{\R^d}\p(t)e^{-2\pi i\la x,y\ra}dt
$$
is a continuous bijection of $\S$ onto $\S$.
The elements of the space $\S'$ of all continuous linear functionals on $\S$ are called {\it tempered} distributions.
For every $\Phi\in\S'$ its Fourier transform
$\hat\Phi$ is defined by the equality
\begin{equation}\label{h}
\hat\Phi(\p)=\Phi(\hat\p)\quad \forall\,\p\in\D.
\end{equation}
Since the space $\D$ is dense in $\S$, we find that $\hat\Phi$ belongs to $\S'$.

A positive measure $\nu$ on $\R^d$ is tempered if and only if 
$$
\log\nu(B(0,r))=O(\log r),\qquad (r\to\infty),
$$ 
where, as usual, $B(x,r)$ denotes the ball with center $x$ and radius $r$ (see, e.g., Lemma 1 from \cite{F1}). It is clear that for a pure point measure
\begin{equation}\label{b}
\hat\mu=\sum_{\g\in\G}b_\g\d_\g
\end{equation}
with a countable $\G$, the measure $|\hat\mu|$ is tempered if and only if
$$
   \log\sum_{\g\in\G,|\g|<r}|b_\g|=O(\log r),\qquad r\to\infty.
$$
A measure $\mu$ in $\R^d$ is {\it translation bounded} if
$$
\sup_{x\in\R^d}|\mu|(B(x,1))<\infty
$$
(see \cite{M1}).  Clearly, every translation bounded measure $\mu$ satisfies the condition
$$
|\mu|(B(x,r))\le Cr^d,\qquad r\ge1,
$$
with some $C<\infty$. In particular, it follows that $|\mu|$ is tempered.

When $\mu$ is a tempered  measure and $\hat\mu$ is defined as in \eqref{b},
the equality \eqref{h} has the form
$$
   \sum_{\g\in\G}b_\g\p(\g)=\int_{\R^d}\hat\p(x)\mu(dx),\qquad \forall\p\in\D.
$$
For $\mu=\sum_{\l\in\L}a_\l\d_\l$ the above equality becomes
$$
   \sum_{\g\in\G}b_\g\p(\g)=\sum_{\l\in\L}a_\l\hat\p(\l),\qquad \forall\p\in\D,
$$
and is called {\it generalized Poisson formula}. If  $a_\l\neq0\ \forall\ \l\in\L$ and $b_\g\neq0\ \forall\g\in\G$,
we say that $\L$ is the support of $\mu$ and $\G$ is the spectrum of $\mu$.
\medskip

With the above notations and properties, we have the following result from \cite{F2}:
\begin{Th}\label{T1} {\rm (\cite{F2})}
Let $\mu$ and $\hat\mu$ be tempered measures on $\R^d$ and $\mu$ have a polynomially discrete support, that is
$$
     |t-t'|\ge c\min\{1,\,|t|^{-h}\}\qquad \forall t, t'\in\supp\mu.
$$
Then $|\mu|$ is tempered. If, in addition, $\hat\mu$ has a polynomially discrete support, then $\mu$ is a Fourier quasicrystal.
\end{Th}

In this paper, we obtain the following results:
\begin{Th}\label{T2}
Let $\mu$ be a positive or translation bounded measure on $\R^d$ and its Fourier transform $\hat\mu$ be a pure point measure \eqref{b}.
Then the measure
\begin{equation}\label{sq}
 \nu=\sum_{\g\in\G}|b_\g|^2\d_\g
\end{equation}
is translation bounded. If $\G$ is locally finite and the numbers $\#(\G\cap B(0,r))$ grow polynomially as $r\to\infty$,
then $|\hat\mu|$ is tempered. If, in addition, $\mu$  is a crystalline measure, then $\mu$ is a Fourier quasicrystal.
\end{Th}

\begin{Th}\label{T3}
Let $\mu$ be a positive or translation bounded measure on $\R^d$ and its Fourier transform $\hat\mu$ be a pure point measure \eqref{b}.
If there exists $\eta>0$ such that the sets $\G\cap B(x,\eta)$ are linearly independent over $\Z$ for all $x\in\R^d$,
then $\hat\mu$ is translation bounded.  In particular, if $\mu$ is a positive crystalline measure
 with a linearly independent over $\Z$ spectrum, then $\mu$ is a Fourier quasicrystal.
\end{Th}
Recall that the set $X\subset\R^d$ is linearly independent over $\Z$ if for any finite number of elements $x_1,\dots,x_n\in X$
and arbitrary $m_1,\dots,m_n\in\Z$ the equality $m_1x_1+\dots+m_nx_n=0$ implies $m_1=\dots=m_n=0$.
Note that crystalline measures with linearly independent over $\Z$ supports were considered by Y. Meyer \cite{M2}.

\begin{Cor}
Let the measure $\mu$ from Theorem \ref{T3} have a uniformly discrete\footnote{A set $A$ is uniformly discrete if
$\inf_{x,x'\in A,x\neq x'}|x-x'|>0.$}
support $\L$ and let $\inf_{\l\in\L}|\mu(\l)|>0$. Then
there exist a positive integer $N$, lattices $L_1,\dots,L_N$ in $\R^d$ of rank $d$ (some of them may coincide), points $\l_1,\dots,\l_N\in\L$,
bounded sets $S_1,\dots,S_N\subset\R^d$  such that the support of $\mu$ is a finite union of shifts of these lattices and
\begin{equation}\label{ad}
\mu=\sum_{j=1}^N\sum_{x\in L_j+\l_j}\left[\sum_s\b_{j,s} e^{2\pi i\la x,\a_{j,s}\ra}\right]\d_x
\end{equation}
with $\a_{j,s}\in S_j$, $\b_{j,s}\in\C$, and $\sum_s|\b_{j,s}|<\infty$ for all $j=1,\ldots,N$.
\end{Cor}

The paper is organized as follows. In Section 2 we collect definitions and properties of almost periodic functions. Three necessary preliminary results are given in Section 3. Section 4 is devoted to the proofs of Theorems \ref{T2} and \ref{T3}.

\bigskip
\section{Almost periodic functions}\label{S2}
\bigskip

The proofs of our theorems are based on some results of the classical theory of almost periodic functions
(see \cite{B} for $d=1$ and \cite{M1}, \cite{R} for arbitrary $d$).

\begin{Def}\label{D1}
A continuous function $f(x)$ on the space $\R^d$
is almost periodic  if for any
  $\e>0$ the set of $\e$-almost periods
  $$
E_{\e}= \{\tau\in\R^d:\,\sup_{x\in\R^d}|f(x+\tau)-f(x)|<\e\}
  $$
is relatively dense, i.e., $E_{\e}\cap B(x,L)\neq\emptyset$ for all $x\in\R^d$ and some $L$ depending on $\e$.
\end{Def}

Each almost periodic function in $\R^d$ is bounded,  a finite sum or product of almost periodic  functions is also almost periodic,
and the uniform limit in $\R^d$ of almost periodic functions is almost periodic too. In particular,
every  absolutely convergent Dirichlet series
\begin{equation}\label{Dir}
   D(x)=\sum_{\o\in\O}a_\o e^{2\pi i\la x,\o\ra},\qquad \sum_{\o\in\O}|a_\o|<\infty,
\end{equation}
with a countable $\O\subset\R^d$ is almost periodic. In this case 
$$
\hat D=\sum_{\o\in\O}a_\o \d_\o.
$$
For every almost periodic $f(x)$ and every $\o\in\R^d$ one considers
the Fourier coefficients $c_\o(f)\in\C$ as 
$$
    c_\o(f)=\lim_{R\to\infty}\frac{1}{v_d R^d}\int_{B(x,R)}f(t)e^{-2\pi i\la t,\o\ra}dt,
$$
where $v_d$ is the volume of the unit ball in $\R^d$. The limit exists uniformly in $x\in\R^d$.
The set $\{\o:\,c_\o(f)\neq 0\}$ is at most countable.
It is easy to check that uniform convergence of the series \eqref{Dir} implies $c_\o(D)=a_\o\ \forall\o\in\R^d$.
Then, an analog of the Parseval identity
\begin{equation}\label{Par}
 \sum_{\o\in\R^d}|c_\o(f)|^2=\lim_{R\to\infty}\frac{1}{v_d R^d}\int_{B(x,R)}|f(t)|^2dt   
\end{equation}
holds. 

Note that by Y.Meyer (cf. \cite[Theorem 3.8]{M3}), if the Fourier transform $\hat f$ of an almost periodic function $f$ is a measure, then it has the form
$$
   \hat f=\sum_{\o\in\O} c_\o(f)\d_\o,\qquad \sum_{|\o|<r}|c_\o(f)|<\infty, \quad \forall\,r<\infty.
$$

\bigskip
\section{Preliminary results}\label{S3}
\bigskip

We shall need the following multidimensional version of the Kroneker lemma (cf. \cite[Chapter III]{Cas}). We provide a proof for self-completeness. 

\begin{Pro}\label{P1}
Let $x_1,\dots,x_N$ be linearly independent over $\Z$ vectors from $\R^d$ and $\t_1,\dots,\t_N$ be arbitrary real numbers. Then for every $\e>0$
there exist $t\in\R^d$ and $p_j\in\Z$, $j=1,\ldots,N$, such that \begin{equation}\label{in}
|\la x_j,t\ra-\t_j-p_j|<\e,\quad j=1,\ldots,N.
\end{equation}
\end{Pro}

{\bf Proof}. We set
$$
 f(t)=1+e^{2\pi i[\la x_1,t\ra-\t_1]}+\dots+e^{2\pi i[\la x_N,t\ra-\t_N]}.
$$
It is enough to prove that $\sup_{t\in\R^d}|f(t)|=N+1$.
For $q\in\N$ we have
\begin{equation}\label{q}
f^q(t)=\sum_s\a_s e^{2\pi i\la\b_s,t\ra},
\end{equation}
where
\begin{equation}\label{p}
\b_s=m_1 x_1+\dots+m_N x_N,
\end{equation}
with $m_j \in\Z$, $0\le m_j\le q$, $j=1,\ldots,N$, and $m_1+\cdots+m_N \leq q$.
Therefore, the number of the terms on the right-hand side of \eqref{q} is at most $(q+1)^N$.

Further, $f^q(t)=(1+z_1+\dots+z_N)^q$, where $z_j=e^{2\pi i[\la x_j,t\ra-\t_j]}$, and every term
 can be written as 
$$
 C_{m_1,\dots,m_N}e^{2\pi i[\la m_1x_1+\dots+m_Nx_N,t\ra-\phi_s]},
$$
where
$$
\phi_s=m_1\t_1+\dots+m_N\t_N,
$$
and the (multi-binomial) coefficient $C_{m_1,\dots,m_N}$ 
coincides with the coefficient of $z_1^{m_1}\cdots z_N^{m_N}$ in the expansion of $(1+z_1+\dots+z_N)^q$ and is, therefore, positive.

It follows from the linear independence of $x_1,\dots,x_N$ that the
monomial corresponding to $\b_s$ is uniquely determined.   
Hence the corresponding coefficient $\a_s$  has the form
\begin{equation}\label{K}
     \a_s=C_{m_1,\dots,m_N}e^{-2\pi i\phi_s},
\end{equation}
 Consequently, we have 
$$
  \sum_s|\a_s|=\sum_{m_1,\dots,m_N}C_{m_1,\dots,m_N}=(N+1)^q.
$$
On the other hand,
$$
\a_s=\lim_{R\to\infty}\frac{1}{v_d R^d}\int_{B(x,R)}e^{-2\pi i\la t,\b\ra}f^q(t)dt.
$$
Therefore, if $\sup_{t\in\R^d}|f(t)|=r<N+1$, then $|\a_s|\le r^q$ and $\sum_s|\a_s|\le r^q(q+1)^N$.
But this is impossible because $r^q(q+1)^N<(N+1)^q$ for large enough $q$. \B
\medskip

{\it Remark 1}. Without the linear independence assumption we may have two distinct representations for some $\b_s$ in \eqref{p}. Then the above argument is valid if a unique $\phi_s$ in \eqref{K} (up to an integer)  corresponds to these two representations (collections of $m_j$). Therefore, we can replace the linear independence condition of $x_j$ with the following: "if $m_1x_1+\dots+m_Nx_N=0$, then $m_1\t_1+\dots+m_N\t_N\in\Z $". It is easy to prove that the latter condition is also necessary for the existence of $t\in\R^d$ in \eqref{in} for all $\e>0$. The one-dimensional version of Kronecker lemma was stated in exactly this way in \cite{Le}.  
\medskip

We also need two results on the connection between measures and their Fourier transforms.

\begin{Pro}\label{P2} {\rm (cf. Theorem 1 from \cite{F6})}.
 If $\mu$ is a nonnegative tempered measure on $\R^d$ and $\hat\mu$ is a complex measure, then $\mu$ is translation bounded.
\end{Pro}

{\bf Proof}.  For a $C^\infty$-function $\psi(t)\not\equiv0$ such that $\supp\psi\in B(0,1)$ we consider 
$$\p(y)=\int_{\mathbb{R}^d}\psi(y-t)\overline{\psi(t)}dt.$$
Then we obtain
$$
\supp\p\subset B(0,2),\quad\hat\p(x)=|\hat\psi(x)|^2\ge0,\quad \hat\p(x)\not\equiv0,
$$
and, therefore, there exist $\eta>0$, $r>0$, and $x_0\in\R^d$, such that
$\hat\p(x)>\eta$ for all $x\in B(x_0,r)$. Then for any $t\in\R^d$ we consecutively have 
\begin{eqnarray*}
  \mu(B(t,r)) &\leq& \eta^{-1}\int_{B(t,r)}\hat\p(x-t+x_0)\mu(dx) \\
  &\leq& \eta^{-1}\int_{\R^d}\hat\p(x-t+x_0)\mu(dx)\\
  &=& \eta^{-1}\int_{\R^d}\p(y)e^{2\pi i(t-x_0)y}\hat\mu(dy) \\ 
  &\leq&  \eta^{-1}\max_{\R^d}|\p(y)||\hat\mu|(B(0,2)).
  \end{eqnarray*}
Since $|\hat\mu|$ is a measure, we conclude that $|\hat\mu|(B(0,2))$ is finite and the measure $\mu$ is translation bounded. \B

\begin{Pro}\label{P3}
 Let $\mu$ be a  translation bounded measure on $\R^d$ and $\p\in\S$. Then
the function $\mu\star\p$ is bounded on $\R^d$ by a constant that depends neither on shifts of $\mu$ and $\p$ nor on their
multiplying by $e^{iat},\ a\in\R$.
\end{Pro}
{\bf Proof}. Set $\mu_x(E)=\mu(E-x)$ for all $E\subset\R^d,\ x\in\R^d$. Then we get
$$
|\mu_x|(B(0,r))\le C_1\max\{1,r\}^d
$$
for all $x\in\R^d$. Since $\p\in\S$, we have $|\p(t)|\le C_2\max\{1,\,r\}^{-d-1}$.

It follows that
\begin{eqnarray*}
|(\mu\star\p)(x)| &=& \left|\int_{\R^d}\p(t)\mu_x(dt)\right| \leq C_2\left[\int_{|t|\le1}|\mu_x|(dt)+\int_{|t|<1}|t|^{-d-1}|\mu_x|(dt)\right] \\
&=& C_2(d+1)\int_1^\infty\frac{|\mu_x|(B(0,s))ds}{s^{d+2}}\le (d+1)C_1 C_2,
\end{eqnarray*}
which completes the proof.
\B

\bigskip
\section{Proofs of Theorems \ref{T2} and \ref{T3}}\label{S4}
\bigskip

In this section, we present proofs of Theorems \ref{T2} and \ref{T3}.

{\bf Proof of Theorem \ref{T2}}. Let $\p\in\S$ be a non-negative even function with compact support such that $\p(y)\equiv1$ for $|y|<1$.
 We have
\begin{eqnarray*} \label{e1t2}
 \mu\star\hat\p(t)&=&\int_{\R^d}\widehat{\left[\p(-y)e^{2\pi i\la t,y\ra}\right]}(x)\mu(dx)\\
 &=& \int_{\R^d}\p(-y)e^{2\pi i\la t,y\ra}\hat\mu(dy) \\
 &=&\sum_{\g\in\supp\p}b_\g\p(\g)e^{2\pi i\la t,\g\ra}.
\end{eqnarray*}
Since $\hat\mu$ is a measure, we see that $\sum_{\g\in\supp\p}|b_\g|$ is finite and the last sum is absolutely convergent. Therefore, $\mu\star\hat\p$ is an almost periodic function. We also obtain from the above equality that 
$$
\widehat{\mu\star\hat\p}=\sum_{\g\in\supp\p}\p(\g)b_\g \d_\g.
$$
By Meyer's result (see the end of Section \ref{S2}), the Fourier coefficients of the function $\mu\star\hat\p$ are
$$
   c_\g(\mu\star\hat\p)=b_\g\p(\g).
$$
Using Parseval's equality \eqref{Par}, we obtain that the variation of the measure $\nu$ from \eqref{sq} can be estimated as follows:
\begin{eqnarray*}
  |\nu|(B(0,1)) &=& \sum_{|\g|<1}|b_\g|^2 \\ \
  &\leq& \sum_{\g\in\R^d}|\p(\g)b_\g|^2\\ &=& \sum_{\g\in\R^d}| c_\g(\mu\star\hat\p)|^2 \\
&=&\lim_{R\to\infty}\frac{1}{v_d R^d}\int_{B(0,R)}|(\mu\star\hat\p)(x)|^2dx.  
\end{eqnarray*}
 It follows from  Proposition \ref{P2}  that under the conditions of the theorem, the measure $\mu$ is translation bounded. Hence by Proposition \ref{P3}, the last integral is bounded by a constant $C$. If we replace $\p(y)$ with $\p(y-y_0)$ (here $y_0\in\R^d$ can be arbitrary), then $\hat\p$ changes by a factor $e^{2\pi i\la x, y_0\ra}$ and by Proposition \ref{P3}, the constant $C$ remains the same. We conclude that
\begin{equation}\label{n}
   \nu(B(y_0,1))\le C.
\end{equation}
Consequently, the measure $\nu$ is translation bounded, and $\nu(B(0,r))\le C'r^d$ with some constant $C'$.

Finally, if
$$
  \#\{\g\in\G:\,|\g|<r\}=O(r^\rho),\qquad (r\to\infty),
$$
with some $\rho<\infty$, then by the Cauchy-Bunyakovskii inequality
\begin{eqnarray*}
     \sum_{|\g|<r}|b_\g| &\leq& \left[\sum_{|\g|<r}|b_\g|^2\right]^{1/2}
     \cdot [\#\{\g\in\G:\,|\g|<r\}]^{1/2} \\
     &=& O(r^{(d+\rho)/2})
\end{eqnarray*}
as $r$ tends to infinity. Therefore, $|\hat\mu|\in\S'$. If $\mu$ is a crystalline measure,
we obtain that it is a Fourier quasicrystal. This completes the proof.  \B
\medskip

{\it Remark 2}. By \eqref{n}, the masses $b_\g$ are uniformly bounded as well. Hence we can replace the exponent $2$ in \eqref{sq} with any $q>2$. We do not know whether \eqref{sq} is valid for exponent $q<2$ under the conditions of Theorem \ref{T2}.

\medskip
{\bf Proof of Theorem \ref{T3}}. It follows from  Proposition \ref{P2}  that under the conditions of the theorem, the measure $\mu$ is translation bounded. Suppose for a contradiction that $\hat\mu$ is not translation bounded. 
Then there exists $\eta>0$ such that the masses of the measure $|\hat\mu|$ in balls of radius $\eta/2$ are unbounded. This means that there are points $y_n\in\R^d$, $n=1,2,\ldots$, such that
$$
    \sum_{\g\in\G\cap B(y_n,\eta/2)}|b_{\g}|>5n.
$$

If the set $\G\cap B(y_n,\eta)$ is finite for some $n\in\N$, then we set $A_n=\G\cap B(y_n,\eta)$. Otherwise, taking into account that
$\hat\mu$ is a measure and $|\hat\mu|(B(y_n,\eta))<\infty$, we
choose $A_n$ to be a finite subset of $\G\cap B(y_n,\eta)$ such that
\begin{equation}\label{sm}
  \sum_{\g\in\G\cap B(y_n,\eta)\setminus A_n}|b_{\g}|<n.
\end{equation}
With this setting, we have
$$
 \sum_{\g\in A_n\cap B(y_n,\eta/2)}|b_{\g}|>4n.
$$
 Applying Proposition \ref{P1}, we can find $x_n \in \R^d$ and numbers $m_\g\in\N$ such that for all $\g\in A_n$
$$
|\la x_n,\g\ra+\frac{\arg b_{\g}}{2\pi}-m_\g|<\frac{\pi}{3}.
$$
Therefore, for all $\g\in A_n$ we have
$$
  \Re e^{2\pi i\la x_n,\g\ra}b_\g> \frac{|b_\g|}{2} >0,
$$
and
\begin{equation}\label{arg}
\Re\sum_{\g\in A_n\cap B(y_n,\eta/2)}e^{2\pi i\la x_n,\g\ra}b_\g>2n.
\end{equation}
Let $\p\in\D$ be such that $0\le\p(y)\le1$, $\supp\p\subset B(0,\eta)$, and $\p(y)\equiv1$ for $|y|<\eta/2$.
It follows from \eqref{sm} and \eqref{arg} that
\begin{multline}\label{c1}
   \left|\int e^{2\pi i\la x_n,y\ra}\p(y-y_n)\hat\mu(dy)\right|    \ge\left|\sum_{\g\in A_n}\p(\g-y_n)e^{2\pi i\la x_n,\g\ra}b_\g\right|-\sum_{\g\in\G\cap B(y_n,\eta)\setminus A_n}|b_{\g}|\\
   \ge \Re\sum_{\g\in A_n\cap B(y_n,\eta/2)}e^{2\pi i\la x_n,\g\ra}b_\g+
   \Re\sum_{\g\in A_n\setminus B(y_n,\eta/2)}\p(\g-y_n)e^{2\pi i\la x_n,\g\ra}b_\g-n\ge n.
\end{multline}
On the other hand,
$$
   \int e^{2\pi i\la x_n,y\ra}\p(y-y_n)\hat\mu(dy)=\int e^{-2\pi i\la y_n,x\ra}\hat\p(x-x_n)\mu(dx)=
   \int e^{-2\pi i\la y_n,x\ra}\hat\p(x)\mu_{x_n}(dx).
$$
By Proposition \ref{P3}, the modulus of the integral in the right-hand side is bounded with a constant that does not depend on
$x_n$ and $y_n$. This inequality contradicts to what is obtained in \eqref{c1}.  Therefore, the measure $\hat\mu$ is translation bounded. If in the conditions of the theorem the measures $\mu$ and $\hat\mu$ have locally finite supports, that is, $\mu$ is a crystalline measure, then $\mu$ is a Fourier quasicrystal.   \B
\medskip

{\bf Proof of Corollary}. It follows from Theorem \ref{T3} that $|\hat\mu|(B(0,r))=O(r^d)$ as $r\to\infty$. Therefore,
all the conditions of Theorem 3 from \cite{F0} are satisfied and we obtain the representation \eqref{ad}.  \B
\bigskip

{\bf Question}. Let $\mu$ be a crystalline measure, and $|\mu|$ be a tempered measure. Is $|\hat\mu|$ a tempered measure?
And what if the measure $\mu$ is translation bounded?

\medskip

{\bf Acknowledgements.} This research was partially supported by Bulgarian Ministry
of Education and Science, Scientific Programme "Enhancing the Research Capacity in Mathematical Sciences (PIKOM)", No. DO1-241/15.08.2023 and by Georgi Chilikov Foundation.

\end{document}